\documentclass[fleqn]{mat01}
\usepackage{times,mathtimy,amssymb}
\begin{document}

\setcounter{page}{217}
\firstpage{217}

\def\theor{\trivlist\item[\hskip\labelsep{{\bf Theorem.}}]}

\newtheorem{theorr}{Theorem}
\renewcommand\thetheorr{\arabic{theorr}}
\newtheorem{therr}[theorr]{\bf Theorem}
\newtheorem{lem}[theorr]{Lemma}
\newtheorem{exampl}[theorr]{Example}
\newtheorem{quest}[theorr]{Question}
\newtheorem{propo}[theorr]{\rm PROPOSITION}
\newtheorem{case}[theorr]{\it Case}

\def\G{G^\prime}
\def\A{X = \{ 1,2,3 \}}
\def\lhd{\mbox{$\, \vartriangleleft\ \ $}}

\title{On finite groups whose every proper normal subgroup is a union of
a given number of conjugacy classes}

\markboth{Ali Reza Ashrafi and Geetha Venkataraman}{On finite groups of
conjugacy classes}

\author{ALI REZA ASHRAFI and GEETHA VENKATARAMAN$^{*}$}

\address{Department of Mathematics, University of Kashan,
Kashan, Iran\\
\noindent $^{*}$Department of Mathematics and Mathematical Sciences
Foundation, St. Stephen's College, Delhi 110 007, India\\
\noindent E-mail: ashrafi@kashanu.ac.ir; geetha\_venkat@ststephens.edu}

\volume{114}

\mon{August}

\parts{3}

\Date{MS received 19 June 2002; revised 26 March 2004}

\begin{abstract}
Let $G$ be a finite group and $A$ be a normal subgroup of $G$. We denote
by $ncc(A)$ the number of $G$-conjugacy classes of $A$ and $A$ is called
$n$-decomposable, if $ncc(A)=n$. Set ${\cal K}_G = \{ncc(A)| A \lhd G
\}$. Let $X$ be a non-empty subset of positive integers. A group $G$ is
called $X$-decomposable, if ${\cal K}_G = X$.

Ashrafi and his co-authors \cite{ash1,ash2,ash3,ash4,ash5} have characterized the
$X$-decomposable non-perfect finite groups for $X = \{ 1, n \}$ and $n
\leq 10$. In this paper, we continue this problem and investigate the
structure of $X$-decomposable non-perfect finite groups, for $X = \{ 1,
2, 3 \}$. We prove that such a group is isomorphic to $Z_6, D_8, Q_8,
S_4$, SmallGroup(20, 3), SmallGroup(24, 3), where SmallGroup$(m,n)$
denotes the $m$th group of order $n$ in the small group library of GAP
\cite{gap}.
\end{abstract}


\maketitle

\section{Introduction and preliminaries}

Let $G$ be a finite group and let ${\cal N}_G$ be the set of proper
normal subgroups of $G$. An element $K$ of ${\cal N}_G$ is said to be
$n$-decomposable if $K$ is a union of $n$ distinct conjugacy classes of
$G$. In this case we denote $n$ by $ncc(K)$. Suppose ${\cal K}_G = \{
ncc(A) | A \in {\cal N}_G \}$ and $X$ is a non-empty subset of
positive integers. A group $G$ is called $X$-decomposable, if ${\cal
K}_G = X$. For simplicity, if $X = \{ 1,n \}$ and $G$ is
$X$-decomposable, then we say that $G$ is $n$-decomposable.

In \cite{shi4}, Wujie Shi defined the notion of a complete normal
subgroup of a finite group, which we call 2-decomposable. He proved that
if $G$ is a group and $N$ a complete normal subgroup of $G$, then $N$ is
a minimal normal subgroup of $G$ and it is an elementary abelian
$p$-group. Moreover, $N\subseteq Z(O_p(G))$, where $O_p(G)$ is a maximal
normal $p$-subgroup of $G$, and $|N| (| N|-1)$ divides
$|G|$ and in particular, $|G|$ is even.

Shi \cite{shi4} proved some deep results about finite group $G$ of
order $p^aq^b$ containing a 2-decomposable normal subgroup $N$. He
proved that for such a group $|N| = 2, 3, 2^{b_1}$ or $2^{a_1}+1$, where
$2^{b_1}-1$ is a Mersenne prime and $2^{a_1}+1$ is a Fermat prime.
Moreover, we have (i) if $|N| = 2$, then $N \subseteq Z(G)$, (ii) if $|N|
= 3$, then $G$ has order $2^a3^b$, (iii) if $|N| = 2^{b_1}$, then $G$ has
order $(2^{b_1}-1)2^b$ and (iv) if $|N| = 2^{a_1}+1$, then $G$ has
order $2^a(2^{a_1}+1)^b$.

Next, Wang Jing \cite{shi5}, continued Wujie Shi's work and defined
the notion of a sub-complete normal subgroup of a group $G$,
which we call 3-decomposable. She proved that if $N$ is a sub-complete
normal subgroup of a finite group $G$, then $N$ is a group in which
every element has prime power order. Moreover, if $N$ is a minimal
normal subgroup of $G$, then $N\subseteq Z(O_p(G))$, where $p$ is a
prime factor of $|G|$. If $N$ is not a minimal normal subgroup of
$G$, then $N$ contains a complete normal subgroup $N_1$, where $N_1$ is
an elementary abelian group with order $p^a$ and we have:

\begin{enumerate}
\renewcommand\labelenumi{(\alph{enumi})}
\leftskip .1pc
\item $N=N_1Q$ has order $p^aq$ and every element of $N$ has prime power
order, $|Q|=q$, $q\not =p, q$ is a prime and $G=MN_1$, $M\cap N_1=1$,
where $M=N_G(Q)$,

\item $N$ is an abelian $p$-group with exponent $\leq p^2$ or a special
group; if $N$ is not elementary abelian, then $N_1\leq \Phi(G)$, where
$\Phi(G)$ denotes the Frattini subgroup of $G$.
\end{enumerate}

Shahryari and Shahabi \cite{shah1,shah2} investigated the structure of
finite groups which contain a 2- or 3-decomposable subgroup. Riese and
Shahabi \cite{shah3} continued this theme by investigating the structure
of finite groups with a 4-decomposable subgroup. Using these works in
\cite{ash1} and \cite{ash2}, Ashrafi and Sahraei characterized the
finite non-perfect $X$-groups, for $X = \{ 1,n \}$, $n \leq 4$. They
also obtained the structure of solvable $n$-decomposable non-perfect
finite groups. Finally, Ashrafi and Zhao \cite{ash3} and Ashrafi and Shi
\cite{ash4,ash5} characterized the finite non-perfect $X$-groups, for $X
= \{ 1,n \}$, where $5 \leq n \leq 10$.

In this paper we continue this problem and characterize the non-perfect
$X$-decomposable finite groups, for $X= \{ 1,2,3 \}$. We prove that such
a group is solvable and determine the structure of these groups. In
fact, we prove the following theorem:

\begin{theor} {\it Let $G$ be a non-perfect $\{ 1,2,3\}$-decomposable
finite group. Then $G$ is isomorphic to $Z_6, D_8, Q_8, S_4,$
${\rm SmallGroup}(20,3)$ or ${\rm SmallGroup}(24,3)$.}\vspace{.6pc}
\end{theor}

Throughout this paper, as usual, $G^\prime$ denotes the derived subgroup
of $G$, $Z_n$ denotes the cyclic group of order $n$, $E(p^n)$ denotes an
elementary abelian $p$-group of order $p^n$, for a prime $p$ and $Z(G)$
is the center of $G$. We denote by $\pi(G)$, the set of all prime
divisors of $|G|$ and $\pi_e(G)$ is the set of all orders of elements of
$G$. A group $G$ is called non-perfect, if $G^\prime \ne G$. Also,
$d(n)$ denotes the set of positive divisors of $n$ and
SmallGroup$(n,i)$ is the $i$th group of order $n$ in the small group library of
GAP \cite{gap}. All groups considered are assumed to be finite. Our
notation is standard and is taken mainly from \cite{conw,hupp,rob,sah}.

\section{Examples}

In this section we present some examples of $X$-decomposable finite
groups and consider some open questions. We begin with the finite
abelian groups.

\begin{lem} Let $G$ be an abelian finite group. Set $X = d(n)-\{ n \},$
where $n = |G|$. Then $G$ is $X$-decomposable.
\end{lem}

\begin{proof}
The proof is straightforward. \hfill $\Box$
\end{proof}

\noindent By the previous lemma a cyclic group of order $n$ is $(d(n) -
\{n\})$-decomposable. In the following examples we investigate the
normal subgroups of some non-abelian finite groups.

\setcounter{theorr}{0}
\begin{exampl}{\rm Suppose that $G$ is a non-abelian group of order $pq$,
in which $p$ and $q$ are primes and $p > q$. It is a well-known fact that
$q | p-1$ and $G$ has exactly one normal subgroup. Suppose that $H$ is
the normal subgroup of $G$. Then $H$ is $(1 +
\frac{p-1}{q})$-decomposable. Set $X = \{ 1, 1 + \frac{p-1}{q} \}$. Then
$G$ is $X$-decomposable.}
\end{exampl}

\begin{exampl}{\rm  Let $D_{2n}$ be the dihedral group of order $2n$, $n
\geq 3$. This group can be presented by
\begin{equation*}
D_{2n} = \langle a,b | a^n = b^2 = 1, b^{-1}ab = a^{-1} \rangle.
\end{equation*}
We first assume that $n$ is odd and $X = \{ \frac{d+1}{2} | d | n \}$.
In this case every proper normal subgroup of $D_{2n}$ is contained in
$\langle a \rangle$ and so $D_{2n}$ is $X$-decomposable. Next we assume
that $n$ is even and $Y = \{ \frac{d+1}{2} | d | n\!: 2\!\!\not| d \} \cup
\{ \frac{d+2}{2} |d| n; 2|d \}$. In this case, we can see that
$D_{2n}$ has exactly two other normal subgroups $H = \langle a^2, b
\rangle$ and $K = \langle a^2, ab \rangle$. To complete the example, we
must compute $ncc(H)$ and $ncc(K)$. Obviously, $ncc(H) = ncc(K)$. If $4
| n$, then $ncc(H) = \frac{n}{4} + 2$ and if \hbox{$4\!\!\not|n$,} then $ncc(H) =
\frac{n+6}{4}$. Set $A = Y \cup \{ \frac{n}{4} + 2 \}$ and $B = Y \cup
\{ \frac{n+6}{4} \}$. Our calculations show that if $4 | n$, then
$D_{2n}$ is $A$-decomposable and if \hbox{$4\!\!\not| n$,} then dihedral group
$D_{2n}$ is $B$-decomposable.}
\end{exampl}

\begin{exampl}{\rm  Let $Q_{4n}$ be the generalized quaternion group of order
$4n$, $n \geq 2$. This group can be presented by
\begin{equation*}
Q_{4n} = \langle a,b | a^{2n} = 1, b^2 = a^n, b^{-1}ab = a^{-1}
\rangle.
\end{equation*}
Set $X = \{ \frac{d+1}{2} | d | n$ and $2\!\!\not | d\} \bigcup \{
\frac{d+2}{2} | d | 2n$ and $2|d \}$ and $Y = X \cup \{ \frac{n+4}{2}
\}$. It is a well-known fact that $Q_{4n}$ has $n+3$ conjugacy classes,
as follows:
\begin{align*}
&\{1\} ; \ \{ a^n \}  ; \ \{ a^r,a^{-r} \}( 1 \leq r \leq n-1);\\[.2pc]
&\{ a^{2j}b | 0 \leq j \leq n-1 \} ; \ \{ a^{2j+1}b | 0
\leq j \leq n-1 \}.
\end{align*}
We consider two separate cases that $n$ is odd or even. If $n$ is odd
then every normal subgroup of $Q_{4n}$ is contained in the cyclic
subgroup $\langle a \rangle$. Thus, in this case $Q_{4n}$ is
$X$-decomposable. If $n$ is even, we have two other normal subgroups
$\langle a^2,b \rangle$ and $\langle a^2, ab \rangle$ which are both
$\frac{n+4}{2}$-decomposable. Therefore, $Q_{4n}$ is $Y$-decomposable.}
\end{exampl}

Now it is natural to generally ask about the set ${\cal K}_G = \{ ncc(A)
| A \lhd G \}$. We end this section with the following question:

\setcounter{theorr}{0}
\begin{quest}{\rm Suppose $X$ is a finite subset of positive integers
containing 1. Is there a finite group $G$ which is $X$-decomposable?}
\end{quest}

\section{Main theorem}

Throughout this section $\A$. The aim of this section is to prove the
main theorem of the paper. We will consider two separate cases in which
$G^\prime$ is 2- or 3-decompos- able. In the following simple lemma, we
classify the $X$-decomposable finite abelian\break groups.

\begin{lem}Let $G$ be an abelian $X$-decomposable finite group. Then $G
\cong Z_6${\rm ,} the cyclic group of order $6$.
\end{lem}

\begin{proof}
Apply Lemma 1. \hfill $\Box$
\end{proof}

\noindent For the sake of completeness, we now define two groups $U$ and $V$ which
we will use later. These groups can be presented by
\begin{align*}
U &= \langle x, y, z | x^3 = y^4 = 1, y^2 = z^2, z^{-1}yz =
y^{-1},\\[.2pc]
&\quad\ x^{-1}yx = y^{-1}z^{-1}, x^{-1}zx = y^{-1}  \rangle,\\[.2pc]
V &= \langle x, y | x^4 = y^5 = 1, x^{-1}yx = y^2 \rangle.
\end{align*}
We can see that $U$ and $V$ are groups of orders 24 and 20 which are
isomorphic to SmallGroup$(24,3)$ and SmallGroup$(20,3)$, respectively.
Also, these groups are $X$-decomposable.

To prove the main result of the paper, we need to determine all of
$X$-decomposable groups of order 8, 12, 18, 20, 24, 36 and 42. The
following GAP program determines all the $X$-decomposable groups of the
mentioned orders.

\vskip 3mm {\bis\tt
\noindent AppendTo("x.txt","Begining the Program","$\backslash$n");\\
E:=[8,12,18,20,24,36,42];\\
for m in E do

\hskip .2cm n:=NrSmallGroups(m);

\hskip .2cm F:=Set([1,2,3]);

\hskip 1.3cm for i in [1,2..n] do

\hskip 1.8cm G1:=[];

\hskip 1.8cm G:=[];

\hskip 1.8cm g:=SmallGroup(m,i);

\hskip 1.8cm h:=NormalSubgroups(g);

\hskip 1.8cm d1:=Size(h);d:=d1-1;

\hskip 2.2cm for j in [1,2..d] do

\hskip 2.7cm s:=FusionConjugacyClasses(h[j],g);

\hskip 2.7cm s1:=Set(s);

\hskip 2.7cm Add(G,s1);

\hskip 2.2cm od;

\hskip 2.2cm for k in G do

\hskip 2.7cm a:=Size(k);

\hskip 2.7cm Add(G1,a);

\hskip 2.2cm od;

\hskip 1.7cm G2:=Set(G1);

\hskip 1.7cm if G2=F then AppendTo("x.txt","S(",m,",",i, ")",

\hskip 1.7cm " ");fi;

\hskip 1.2cm od;

\noindent od;
}

\setcounter{theorr}{0}
\begin{propo}$\left.\right.$\vspace{.5pc}

\noindent Let $G$ be a non-perfect and non-abelian $X$-decomposable
finite group such that $\G$ is $2$-decomposable. Then $G$ is isomorphic
to $D_8, Q_8$ or $SmallGroup(20,3)$.
\end{propo}

\begin{proof} Set $G^\prime = 1 \cup Cl_G(a)$. Then it is an easy fact
that $\G$ is an elementary abelian $r$-subgroup of $G$, for a prime $r$.
First of all, we assume that $|\G| = 2$. Then one can see that $\G =
Z(G)$. If $G$ is not a 2-group then there exists an element $x \in G$ of
an odd prime order $q$. Suppose $H = \G\langle x \rangle$. Since $H$ is
a cyclic group of order $2q$, $ncc(H) \geq 4$ which is impossible. Hence
$G$ is a 2-group. We show that $|G| = 8$. Suppose $|G| \geq 16$. Since
$|\G| = 2$ and every subgroup containing $G^\prime$ is normal, we can
find a chain $G^\prime < H < K < G$ of normal subgroups of $G$, a
contradiction. So $G \cong D_8$ or $Q_8$ and by Examples 2 and 3, these
groups are $X$-decomposable.

We next assume that $|\G| \geq 3$. If $Z(G) \neq 1$ then it is easy to
see that $|Z(G)| = 2$ or 3. Suppose $|Z(G)| = 3$. Then $\G{Z(G)} = G$ or
$\G$. If $\G{Z(G)} = G$, then $G \cong \G \times Z(G)$. This implies that
$G$ is abelian, a contradiction. Thus $Z(G) \leq \G$. This leads to a
contradiction, since $G$ is non-perfect and $\G$ is 2-decomposable. Thus
$|Z(G)| = 2$. But in this case $H = Z(G)\G$ is 3-decomposable, which is
impossible. Therefore, $Z(G) = 1$ and by Theorem 2.1 of \cite{shah1}, we
have $|G| = |\G|(|\G|-1)$ and $G$ is a Frobenius group with kernel $\G$
and its complement is abelian. Suppose $|\G| = r^n$, then $|G| =
r^n(r^n-1)$. Take $K$ to be any proper non-trivial subgroup of $T$,
where $T$ is a Frobenius complement of $\G$. Then $K\G$ is
3-decomposable and for any $x \ne 1$ in $K$ we have $|Cl_G(x)| = |G|/|T|
= |\G|$. So $|K\G| = 2|\G|$ and therefore we get $|K| = 2$. Since $T$ is
abelian, this forces 2 to be the only proper divisor of $|T|$ and hence
$|T| = 4$. So $|G| = 20$ and clearly $G$ is a semidirect product of
$Z_5$ by $T$. Further, $Z(G) = 1$ forces $T \leq {\rm Aut}(Z_5) \cong
Z_4$. Hence $T = {\rm Aut}(Z_5)$. Therefore $G \cong {\rm Aut}(Z_5)
\propto Z_5$. To complete the proof, we show that $G \cong $
SmallGroup(20, 3) and it is $X$-decomposable. Let $x$ and $y$ be
elements of $G$ with $o(x) = 4$ and $o(y) = 5$. Since $G$ is a
centerless group containing five involutions, it has exactly two
non-trivial, proper normal subgroups $A = \langle y \rangle$ and $B =
A\langle x^2 \rangle$ of orders 5 and 10, respectively. Clearly $B$ is
non-abelian and so it is isomorphic to the dihedral group of order 10.
This shows that the elements of $B-A$ are conjugate in $B$. But
$x^{-1}yx = y^i$, $i = 2,3,4$. Suppose $i=4$. Since $B \cong D_{10}$, we
have that $x^2yx^{-2} = y^{-1}$. Thus we get that $xyx^{-1} =
x^{-1}y^{-1}x$. Consequently if $i=4$ then we have that $x^{-1}y^{-1}x =
y^{-1}$ and then $G$ would be abelian, a contradiction. Also the two
groups constructed by $i=2$ and $i=3$ will be isomorphic. Hence without
loss of generality, we can assume that $i=2$ and so $G \cong V \cong $
SmallGroup(20, 3). This shows that non-identity elements of $A$ will be
conjugate in $G$ and so $A$ is 2-decomposable and $B$ is
3-decomposable. This completes the proof. \hfill $\Box$
\end{proof}

\begin{propo}$\left.\right.$\vspace{.5pc}

\noindent Let $G$ be a non-perfect and non-abelian $X$-decomposable
finite group such that $\G$ is $3$-decomposable. Then $G$ is isomorphic
to $S_4$ or ${\rm SmallGroup}(24,3)$.
\end{propo}

\begin{proof} Set $G^\prime = 1 \cup Cl_G(a) \cup Cl_G(b)$. Our main
proof will consider three separate\break cases.

\setcounter{theorr}{0}
\begin{case}{\rm $a^{-1} \not\in Cl_G(a)$.\ \ In this case $Cl_G(b) =
Cl_G(a^{-1})$ and by Proposition 1 of \cite{shah2}, $\G$ is an
elementary abelian $p$-subgroup of $G$, for some odd prime $p$. Suppose
$H$ is a 2-decomposable subgroup of $G$. Then by Corollary 1.7 of
\cite{shah1}, we can see that $H = Z(G)$ has order two. Thus $G \cong
Z(G) \times \G$, which is a contradiction.}
\end{case}

\begin{case}{\rm $a^{-1} \in Cl_G(a), b^{-1} \in Cl_G(b)$ and $(o(a),o(b)) =
1$.\ \ In this case, by (\cite{shah2}, Lemma 6), $|\G| = pq^n$, for some
distinct primes $p, q$, and by Lemma 4 of \cite{shah2}, $Z(\G) = 1$. Also,
by  Lemma 5 of \cite{shah2}, $G^{\prime\prime} = 1 \cup Cl_G(a)$ has order
$q^n$. Since $\G$ is 3-decomposable, $|G:\G| = r$, $r$ is prime. Thus
$|G| = prq^n$ and $|\pi(G)| = 2$ or 3. Suppose $|\pi(G)| = 2$. Then by
Shi's result \cite{shi4}, mentioned in the introduction, $|G^{\prime\prime}|
= 2, 3, 2^{b_1}$ or $2^{a_1}+1$, where $2^{b_1}-1$ is a Mersenne prime
and $2^{a_1}+1$ is a Fermat prime. If $|G^{\prime\prime}| = 2$, then $\G$
is a cyclic group of order $2p$, a contradiction. If $|G^{\prime\prime}|
= 3$, then $\G$ is a cyclic group of order $3p$ or isomorphic to the
symmetric group on three symbols. Since $\G$ is centerless, $\G \cong
S_3$. This shows that $|G| = 12$ or $18$ and by our program in GAP
language, there is no $X$-decomposable group of order 12 or 18. We now
assume that $|G^{\prime\prime}| = 2^{b_1}$. Hence by Shi's result,
mentioned before, $|G| = 2^b(2^{b_1}-1)$ and $|\G| = 2^{b_1}(2^{b_1} -
1)$. Suppose $x \in \G - G^{\prime\prime}$ and $y \in G^{\prime\prime}$.
Then we can see that $|Cl_G(x)| = 2^{b_1}(2^{b_1} - 2)$ and $|Cl_G(y)| =
2^{b_1} - 1$. Since $|Cl_G(x)|$ is a divisor of $|G|$, $2^{b_1-1} - 1 |
2^{b_1} - 1$, which implies that $b_1 = 2$. Therefore $|\G| = 12$ and
$|G| = 24$ or 36. Again using our GAP program, we can see that $|G| =
24$ and since $\G$ is centerless, $G \cong S_4$. Next we suppose that
$|G^{\prime\prime}| = 2^{a_1}+1$. Apply Shi's result again to obtain
$|G| = 2^a(2^{a_1}+1)^b$. So $\G$ is a dihedral group of order
$2(2^{a_1}+1)$ and $|G|= 2(2^{a_1}+1)^2$ or $4(2^{a_1}+1)$. If $|G|=
2(2^{a_1}+1)^2$ then $G$ has a 3-decomposable subgroup of order
$(2^{a_1}+1)^2$. This subgroup has a $G$-conjugacy class of length $q^2
- q = 2^{a_1}(2^{a_1}+1)$ and so $a_1 = 1$ and $|G| = 18$, a
contradiction. If $|G| = 4(2^{a_1}+1)$, then $q(q-1) | 4(2^{a_1}+1)$ and
so $a_1 = 1, 2$. This implies that $|G| = 12, 20$, which is impossible.

Therefore it is enough to assume that $|G| = prq^n$, for distinct primes
$p, q$ and $r$. Since $G^{\prime\prime}$ is a 2-decomposable subgroup of
order $q^n$, $q^n(q^n-1)| |G|$. Thus $q^n-1 | pr$. Suppose $n=1$. Then
$|G| = pqr$, $|\G| = pq$ and $\G$ has two $G$-conjugacy classes of
lengths $q-1$ and $q(p-1)$. Hence $p-1 | r$ and $q-1 | pr$. If $p=2$,
then $G$ has a 3-decomposable subgroup $H$ of order $qr$. Since $H$ has
a $G$-conjugacy class of length $q(r-1)$, $r=3$. Thus $q-1 | 6$. This
shows that $q=7$ and $|G| = 42$. But by our GAP program, there is no
$X$-decomposable group of order 42, a contradiction. If $p \ne 2$, then
$p=3, r=2$ and a similar argument shows that $|G| = 42$, which is
impossible. Thus $n \ne 1$. We now assume that $q \ne 2$. Since $q^n-1 |
pr$, $q-1 = p$ or $r$. This shows that $q=3$ and one of $p$ or $r$ is
equal to 2. If $p$ or $r$ take the value 2, using arguments similar to
the case $n=1$, we get $r$ or $p$ equals 3 respectively. This is a
contradiction as $q=3$. Finally, we assume that $q=2$. Then $|\G| =
2^np$ and $|G| = 2^npr$. Since $\G$ has a $G$-conjugacy class of length
$2^n(p-1)$, $p-1 | r$. Therefore, $p=2$ or $r=2$ which is our final
contradic-\break tion.}\vspace{-.2pc}
\end{case}

\begin{case}{\rm $a^{-1} \in Cl_G(a), b^{-1} \in Cl_G(b)$ and $(o(a),
o(b)) \neq 1$.\ \ In this case by Proposition 2 of \cite{shah2}, we have
that $\G$ is a metabelian $p$-group. Since $\G$ is a maximal subgroup
of $G$, we have that $|G:\G| = q$, where $q$ is prime. If $q=p$ then $G$
is $p$-group and so $\G \leq \Phi(G)$. This shows that $G$ is cyclic, a
contradiction. Thus $|G| = p^nq$, for distinct primes $p$ and $q$.
Suppose $H$ is a 2-decomposable subgroup of $G$. If $H$ is central then
$H = Z(G)$. We first assume that $Z(G) \not\leq \G$. Then $G \cong \G
\times Z(G)$ which implies that $Z(\G) = 1$, a contradiction. Next,
suppose that $Z(G) \leq \G$. Then $\G$ has a $G$-conjugacy class of
length $2^n-2$ and so $p=2$ and $q=2^{n-1}-1$. Without loss of
generality we can assume that $|Cl_G(a)| = 1$ and $|Cl_G(b)| = 2^n - 2$.
This shows that $|C_G(b)| = |G|/|Cl_G(b)| = 2^{n-1}$ and $G^\prime$
cannot be abelian because if $\G$ is abelian we would get $\G \leq
C_G(b)$. Since $\G$ is non-abelian, $o(b) = 4$ and $\G$ has a unique
subgroup of order 2. Thus $\G \cong Q_8$ and $G$ is a semidirect product
of $Q_8$ by $Z_3$. Assume that $H = \langle x \rangle$ is the cyclic
group of order 3 and $N = Q_8$. Then ${\rm Aut}(N)$ has a unique
conjugacy type of automorphism of order 3. Therefore $G = H
\propto_\theta N$, where $\theta(x)$ is an automorphism of order 3 of
the group $Q_8$. We now show that $G$ is $X$-decomposable and it is
isomorphic to SmallGroup$(24,3)$. The possible orders for a non-trivial
proper normal subgroups of $G$ are 2, 4, 6, 8, 12. Clearly $Z(G) = Z(N)$
and $G/Z(G) \cong A_4$. But $A_4$ does not have normal subgroups of
order 2, 3 and 6, so every normal subgroup of $G$ has order 2 or 8 and
these are unique. On the other hand, using Example 3, we can assume that
$N = \langle y,z| y^4=1, y^2=z^2, z^{-1}yz=y^{-1} \rangle$. Thus
$G$ is isomorphic to a group which has the same presentation as the
group $V$, which we defined before. This shows that $G \cong {\rm
SmallGroup}(24,3)$ and also $N$ is 3-decompos-\break able.

Now it remains to investigate the case that $H$ is not central. Thus $H
\leq \G$ and $\G$ has a $G$-conjugacy class of length $p^i-1 $, for some
$1 \leq i \leq n-1$. Thus $p^i -1 = q$. If $p=3$, then $i=1, q=2$ and we
can show that $|G| = 18$. Thus $|\G| = 9$ and so $\G$ is abelian.
Further $|H| = 3$, so assuming without loss of generality that $b \in \G
\backslash H$, we get $|Cl_G(b)| = 6$. Hence $|C_G(b)| = 3$ which is a
contradiction as it must be at least 9. Hence $p=2$, $i$ is prime and
$|G| = 2^{2i}(2^i-1)$ or $2^{1+i}(2^i-1)$. Suppose $|G| =
2^{1+i}(2^i-1)$, $Q$ is a Sylow $q$-subgroup of $G$ and $N=HQ$. Since
$G/N$ is abelian, $\G \leq N$, which is impossible as $|\G|$ does not
divide $|N|$. Finally, we assume that $|G| = 2^{2i}(2^i-1)$. We may
assume without loss of generality that $H = 1 \cup Cl_G(a)$. Then we get
that $|Cl_G(a)| = 2^i-1$ and so $|C_G(a)| = 2^{2i}$. Thus $C_G(a) = \G$.
Hence we get $H \leq Z(\G)$ and so $H \leq C_G(b)$. But $|C_G(b)| = 2^i
= |H|$ and so $H = C_G(b)$. Thus $b \in H$, which is our final
contradiction. This completes the proof.}\hfill $\Box$
\end{case}\vspace{-1pc}
\end{proof}

Now we are ready to prove our main result.

\begin{theor}{\it Let $G$ be a non-perfect $X$-decomposable
finite group. Then $G$ is isomorphic to $Z_6, D_8, Q_8, S_4,
{\rm SmallGroup}(20,3)$ or ${\rm SmallGroup}(24,3)$.}
\end{theor}

\begin{proof} The proof is straightforward and follows from Lemma 2,
Proposition 1 and Proposition 2. \hfill $\Box$
\end{proof}

We end this paper with the following question:

\setcounter{theorr}{1}
\begin{quest}{\rm Is there any classification of perfect
$X$-decomposable finite groups?}
\end{quest}

\section*{Acknowledgement}

One of the authors (ARA) expresses his gratitude to Prof Wujie Shi for
pointing out refs \cite{shi4} and \cite{shi5} and is also
obliged to Prof Zhao Yaoqing for the translation of the main results of
\cite{shi4} and \cite{shi5} from Chinese to English.

\end{document}